\documentclass[12pt]{article}
\usepackage{hyperref}
\usepackage{amsfonts}
\usepackage{enumerate}
\usepackage{graphicx}
\begin{document}

\begin{center}
{\large\bf On a boundary value problem of a singular discrete time system with a singular pencil}

\vskip.20in

Ioannis K. Dassios$^{1,2}$\\[2mm]
{\footnotesize
$^{1}$MACSI, Department of Mathematics \& Statistics, University of Limerick, Ireland\\[5pt]
$^{2}$Electricity Research Centre, University College Dublin, Ireland}
\end{center}

{\footnotesize
\noindent
\textbf{Abstract.} In this article, we study a boundary value problem of a class of singular linear discrete time systems whose coefficients are non-square constant matrices or square with a matrix pencil which has an identically zero determinant. By taking into consideration that the relevant pencil is singular, we provide necessary and sufficient conditions for existence and uniqueness of solutions. In addition, a formula is provided for the case of unique solutions and optimal solutions are studied for the cases of no solutions and infinite many solutions. Finally, based on a singular discrete time real dynamical system, numerical examples are given to justify our theory.\\
\\[3pt]
\\[3pt]

\vskip.2in
\section{Introduction}
Linear discrete time systems (or linear matrix difference equations), are systems in which the variables take their value at instantaneous time points. Discrete time systems differ from continuous time ones in that their signals are in the form of sampled data. With the development of the digital computer, the discrete time system theory plays an important role in control theory. In real systems, the discrete time system often appears when it is the result of sampling the continuous-time system or when only discrete data are available for use. Discrete time systems have many applications in economics, physics, circuit theory and other areas. For example in finance, there is the very famous Leondief model, see \cite{2}, or the very important Leslie population growth model and backward population projection, see also \cite{2}. In physics the Host-parasitoid Models, see \cite{11}. Applications of absorbing Markov chains or the distribution of heat through a long rod or bar are other interesting applications suggested in \cite{26}. Thus many authors have studied discrete time systems and their applications, see \cite{1}, \cite{2}, \cite{3}, \cite{4}, \cite{5}, \cite{6}, \cite{7}, \cite{8}, \cite{14}, \cite{15}, \cite{19}, \cite{10}, \cite{12}, \cite{13}, \cite{16}, \cite{17}, \cite{18}, \cite{20}, \cite{21}, \cite{22}, \cite{24}, \cite{25}, \cite{26}, \cite{27}, \cite{28}, \cite{29}.

In this article we study a boundary value problem (BVP) of a class of singular linear matrix difference equations whose coefficients are either non-square constant matrices or square with a matrix pencil which has an identically zero determinant. We consider the singular system
\begin{equation}
FY_{k+1}=GY_k,\quad k=k_0,...,k_N,
\end{equation}
with known boundary conditions of type
\begin{equation}
AY_{k_0}+BY_{k_N}=D.
\end{equation}
Where $F,G \in \mathbb{C}^{r \times m}$, $A,B \in \mathbb{C}^{n \times m}$, $Y_k\in \mathbb{C}^{m \times 1}$ and $D\in \mathbb{C}^{n \times 1}$. System (1) is singular, i.e. the matrices $F$ and $G$ can be non-square ($r\neq m$) or square ($r = m$) and $F$ singular (det$F$=0). Furthermore, $k$ belongs to the set $\left\{k_0,k_0+1,k_0+2,...,k_N-1,k_N,k_N+1,...\right\}\subseteq \mathbb{N}$.

Many authors use matrix pencil theory to study linear discrete time systems with constant matrices, see for instance \cite{4}, \cite{5}, \cite{6}, \cite{7}, \cite{8}, \cite{14}, \cite{15}, \cite{19}, \cite{10}, \cite{13}, \cite{16}.  A matrix pencil is a family of matrices $sF-G$, parametrized by a complex number $s$, see \cite{9}, \cite{13}, \cite{16}, \cite{23}, \cite{25}. When $G$ is square and $F=I_m$, where $I_m$ is the identity matrix, the zeros of the function det($sF-G$) are the eigenvalues of $G$. Consequently, the problem of finding the non-trivial solutions of the equation
\[
sFX=GX
\]
is called the generalized eigenvalue problem. Although the generalized eigenvalue problem looks like a simple generalization of the usual eigenvalue problem it exhibits some important differences. First it is possible for $F$ to be singular in which case the problem has infinite eigenvalues. To see this write the generalized eigenvalue problem in the reciprocal form
\[
FX=s^{-1}GX.
\]
If $F$ is singular with a null vector $X$, then $FX=0_{m, 1}$, so that $X$ is an eigenvector of the reciprocal problem corresponding to eigenvalue $s^{-1}=0$; i.e., $s=\infty$. A second non-trivial case is the determinant det($sF-G$), when $F$, $G$ are square matrices, to be identically zero, independent of $s$. And finally there is the case for both matrices $F$, $G$ to be non-square (for $r\neq m$).

In this article we will consider these last two cases. Actually, we generalize various results regarded the literature which mainly are dealing with square and non-singular systems and apply to general non-square pencils. Another important characteristic of the singular case considered here is that existence of solutions is not automatically satisfied. Explicit and easily testable conditions are derived for which the system has a unique solution. This is very important for many applications for which the model is significant only for certain range of its parameters. In these cases a careful interpretation of results or even a redesign of the system maybe needed. 

The paper is organized as follows: In Section 2 we refer to the mathematical background used throughout this paper, Section 3 contains the main results of the article and Section 4 the numerical examples.

\section{Mathematical background}

In this section we will give the mathematical background and the notation that is used throughout the paper
\\\\
\textbf{Definition 2.1.} Given $F,G\in \mathbb{C}^{r\times m}$ and an arbitrary $s\in\mathbb{C}$, the matrix pencil $sF-G$ is called:
\begin{enumerate}
\item Regular when  $r=m$ and  det$(sF-G)\neq 0$;
\item Singular when  $r\neq m$ or  $r=m$ and det$(sF-G)\equiv 0$.
\end{enumerate}
In this article, we consider the case that the pencil is \emph{singular}.
\\\\
The class of $sF-G$ is characterized by a uniquely
defined element, known as complex Kronecker canonical form,
$sF_K -G_K$, see \cite{9}, \cite{13}, \cite{16}, \cite{23}, specified by the complete set of
invariants of the singular matrix pencil $sF-G$. These invariants are the \emph{elementary divisors} (e.d.) and the \emph{minimal indices} (m.i.). The set of e.d. is obtained by
factorizing the invariant polynomials into powers of homogeneous polynomials irreducible over $\mathbb{C}$. There are two different types of e.d., the set of type  $(s-a_j)^{p_j}$, called \emph{finite elementary divisors} (f.e.d.), where $a_j$ finite eigenvalue of algebraic multiplicity $p_j$ ($1\leq j \leq \nu$) with $\sum_{j =1}^\nu p_j  = p$
, and the set of type $\hat{s}^q=\frac{1}{s^q}$, called \emph{infinite elementary divisors} (i.e.d.), where $q$ is the algebraic multiplicity of the infinite eigenvalue. The set of the m.i. is defined as follows. The distinguishing feature of a singular pencil $sF-G$ is that either $r\neq m$ or $r=m$ and det$(sF-G)\equiv 0$. Let $\mathcal{N}_r$, $\mathcal{N}_l$ be right, left null space of a matrix respectively. Then the equations
\[
(sF-G)U(s)=0_{r, 1}
\]
and
\[
V^T(s)(sF-G)=0_{1, r}.
\]
Where $()^T$ is the transpose tensor, have solutions in $U(s)$, $V(s)$, which are vectors in the rational vector spaces $\mathcal{N}_r(sF-G)$ and $\mathcal{N}_l(sF-G)$ respectively. The binary vectors $U(s)$ and $V^T(s)$ express dependence relationships among the columns or rows of $sF-G$ respectively. $U(s)$, $V(s)$ are polynomial vectors. Let $d$=dim$\mathcal{N}_r(sF-G)$ and $t$=dim$\mathcal{N}_l(sF-G)$. It is known, see \cite{9}, \cite{13}, \cite{16}, \cite{23}, that $\mathcal{N}_r(sF-G)$ and $\mathcal{N}_l(sF-G)$ as rational vector spaces are spanned by minimal polynomial bases of minimal degrees $\epsilon_1=\epsilon_2=...=\epsilon_g=0<\epsilon_{g+1}\leq...\leq\epsilon_d$ and $\zeta_1=\zeta_2=...=\zeta_h=0<\zeta_{h+1}\leq...\leq\zeta_t$ respectively. The set of minimal indices $\epsilon_1, \epsilon_2,..., \epsilon_d$ and $\zeta_1, \zeta_2,..., \zeta_t$ are known as \textit{column minimal indices} (c.m.i.) and \textit{row minimal indices} (r.m.i.) of $sF-G$ respectively. To sum up, in the case of a singular matrix pencil we have invariants of the following type:
\begin{itemize}
    \item e.d. of type  $(s-a_j)^{p_j}$, called \emph{finite elementary
    divisors};
    \item e.d. of type  $\hat{s}^q=\frac{1}{s^q}$, called \emph{infinite elementary divisors};
    \item m.i. of type $\epsilon_1=\epsilon_2=...=\epsilon_g=0<\epsilon_{g+1}\leq...\leq\epsilon_d$, called \emph{column minimal indices};
    \item m.i. of type $\zeta_1=\zeta_2=...=\zeta_h=0<\zeta_{h+1}\leq...\leq\zeta_t$, called \emph{row minimal indices}.
\end{itemize}
\textbf{Definition 2.2.} The direct sum
denoted by $B_{n_1} \oplus B_{n_2}  \oplus \dots \oplus B_{n_r}$ is the block diagonal matrix $blockdiag\left[\begin{array}{cccc} B_{n_1}& B_{n_1}& \dots& B_{n_r}\end{array}\right]$. Where $B_{n_1}\in\mathbb{C}^{n_1\times n_1}$ ,$B_{n_2} \in \mathbb{C}^{n_2\times n_2}$, $\dots$, $B_{n_r}\in\mathbb{C}^{n_r\times n_r}$. 
\\\\ 
The existence of a complete set of invariants for singular pencils implies the existence of canonical form, known as Kronecker canonical form, see \cite{9}, \cite{13}, \cite{16}, \cite{23}, defined by 
\begin{equation}
sF_K  - G_K :=
sI_p  - J_p  \oplus sH_q  - I_q \oplus sF_{\epsilon}-G_{\epsilon}\oplus sF_{\zeta}-G_{\zeta}\oplus 0_{h, g},
\end{equation}
where $p+q+\epsilon+\zeta+h=r$ and $p+q+\epsilon+\zeta+g=m$. The block $sI_p  - J_p$ is uniquely defined by the set of f.e.d. of $sF-G$  
\[
  ({s - a_1 })^{p_1 } , \dots ,({s - a_\nu  }
 )^{p_\nu },\quad \sum_{j = 1}^\nu  {p_j  = p}
\]
and has the form
\[
    sI_p  - J_p  := sI_{p_1 }  - J_{p_1 } (
    {a_1 }) \oplus  \dots  \oplus sI_{p_\nu  }  - J_{p_\nu  }({a_\nu  }).
\]
The $q$  blocks of the second uniquely defined block
$sH_q -I_q$ correspond to the i.e.d. of $sF-G$  
\[
  \hat s^{q_1} , \dots ,\hat s^{q_\sigma}, \quad \sum_{j =
  1}^\sigma  {q_j  = q}
\]
and has the form
\[
    sH_q  - I_q  := sH_{q_1 }  - I_{q_1 }  \oplus
    \dots  \oplus sH_{q_\sigma  }  - I_{q_\sigma}.
\]
The matrix $H_q$  is a nilpotent element of $\mathbb{C}^{q\times q}$  with index $q_* = \max \{ {q_j :j = 1,2, \ldots ,\sigma } \}$, i.e. 
\[
    H^{q_*}_q=0_{q, q}.
\]
In the above notations, the matrices $I_{p_j } ,J_{p_j } ({a_j }),H_{q_j }$ are defined by
\[
   I_{p_j }  = \left[\begin{array}{ccccc} 
   1&0& \ldots & 0&0\\
   0& 1 &  \ldots&0 &0 \\
   \vdots & \vdots & \ddots & \vdots &\vdots \\
   0 & 0 & \ldots  & 0 &1
   \end{array}\right]
   \in \mathbb{C}^{p_j\times p_j } , 
  \]
	\[
   J_{p_j } ({a_j }) =  \left[\begin{array}{ccccc}
   a_j  & 1 & \dots&0  & 0  \\
   0 & a_j  &   \dots&0  & 0  \\
    \vdots  &  \vdots  &  \ddots  &  \vdots  &  \vdots   \\
   0 & 0 &  \ldots& a_j& 1\\
   0 & 0 & \ldots& 0& a_j
   \end{array}\right] \in \mathbb{C}^{p_j\times p_j },
\]
\[
 H_{q_j }  = \left[
\begin{array}{ccccc} 0&1&\ldots&0&0\\0&0&\ldots&0&0\\\vdots&\vdots&\ddots&\vdots&\vdots\\0&0&\ldots&0&1\\0&0&\ldots&0&0
\end{array}
\right] \in \mathbb{C}^{q_j\times q_j }.
  \]
For algorithms about the computations of the Jordan matrices see \cite{9}, \cite{13}, \cite{16}, \cite{23}, \cite{25}. For the rest of the diagonal blocks of (3), i.e. the blocks $sF_{\epsilon}-G_{\epsilon}$ and $sF_{\zeta}-G_{\zeta}$, are defined as follows. The matrices $F_{\epsilon}$, $G_{\epsilon}$ have the form
\begin{equation} 
F_\epsilon=blockdiag\left\{L_{\epsilon_{g+1}}, L_{\epsilon_{g+2}}, ..., L_{\epsilon_d}\right\}
\end{equation}
and
\begin{equation}
G_\epsilon=blockdiag\left\{\bar L_{\epsilon_{g+1}}, \bar L_{ \epsilon_{g+2}}, ..., \bar L_{\epsilon_d}\right\},
\end{equation}  
where $L_{\epsilon_i}= \left[
\begin{array}{ccc} I_{\epsilon_i} & \vdots & 0_{{\epsilon_i}, 1}
\end{array}
\right]$, and $\bar L_{\epsilon_i}=\left[
\begin{array}{ccc} 0_{{\epsilon_i}, 1} & \vdots & I_{\epsilon_i}
\end{array}
\right]$, for $i=g+1, g+2,..., d$. The matrices $F_{\zeta}$, $G_{\zeta}$ have the form
\begin{equation} 
F_\zeta=blockdiag\left\{L_{\zeta_{h+1}}, L_{\zeta_{h+2}}, ..., L_{\zeta_t}\right\}.
\end{equation}
 and
\begin{equation}
G_\zeta=blockdiag\left\{\bar L_{\zeta_{h+1}}, \bar L_{ \zeta_{h+2}}, ..., \bar L_{\zeta_t}\right\},
\end{equation}
where $L_{\zeta_j}= \left[
\begin{array}{c} I_{\zeta_j} \\ 0_{1, \zeta_j}
\end{array}
\right]$, $\bar L_{\zeta_j}=\left[
\begin{array}{c} 0_{1, \zeta_j}\\I_{\zeta_j}
\end{array}
\right]$, for $j=h+1, h+2,..., t$. For algorithms about the computations of these matrices see \cite{9}, \cite{13}, \cite{16}, \cite{23}.

\section{Main results}

Following the analysis in Section 2, there exist non-singular matrices $P$, $Q$ with $P\in \mathbb{C}^{r\times r}$, $Q\in \mathbb{C}^{m\times m}$, such that 
\begin{equation}
PFQ=F_K
\end{equation}
and
\begin{equation}
PGQ=G_K,
\end{equation}
where $F_K$, $G_K$, are defined in (3). Let 
\begin{equation}
Q=\left[\begin{array}{ccccc}Q_p & Q_q &Q_\epsilon & Q_\zeta & Q_g\end{array}\right],
\end{equation}
where $Q_p\in \mathbb{C}^{m\times p}$, $Q_q\in \mathbb{C}^{m\times q}$, $Q_\epsilon\in \mathbb{C}^{m\times \epsilon}$, $Q_\zeta\in \mathbb{C}^{m\times \zeta}$ and $Q_g\in \mathbb{C}^{m\times g}$.
\\\\
\textbf{Lemma 3.1.}
System (1) is divided into five subsystems. The subsystem
\begin{equation}
    Z_{k+1}^p = J_p Z_k^p,
\end{equation}
the subsystem
\begin{equation}
    H_q Z^q_{k+1} = Z_k^q,
\end{equation}
the subsystem
\begin{equation}
    F_\epsilon Z^\epsilon_{k+1}=G_\epsilon Z^\epsilon_k,
\end{equation}
the subsystem
\begin{equation}
    F_\zeta Z^\zeta_{k+1}=G_\zeta Z^\zeta_k,
\end{equation}
and the subsystem
\begin{equation}
    0_{h, g}\cdot Z^g_{k+1}=0_{h, g}\cdot Z^g_k.
\end{equation}
\textbf{Proof.} Consider the transformation
\begin{equation}
    Y_k=QZ_k.
\end{equation}
Substituting (16) into (1) we obtain
\[
    FQZ_{k+1}=GQZ_k,
\]
whereby multiplying by $P$ and using (8) and (9) we arrive at
\[
    F_KZ_{k+1}=G_K Z_k.
\]
Moreover, we can write $Z_k$ as
\[
Z_k=\left[
\begin{array}{c} Z_k^p\\Z_k^q\\Z_k^\epsilon\\Z_k^\zeta\\Z_k^g
\end{array}
\right],
\]
where $Z_k^p\in \mathbb{C}^{p\times 1}$, $Z_k^q \in \mathbb{C}^{q\times 1}$, $Z_k^\epsilon\in \mathbb{C}^{\epsilon\times 1}$, $Z_k^\zeta \in \mathbb{C}^{\zeta\times 1}$ and $Z_k^g\in \mathbb{C}^{g\times 1}$. By taking into account the above expressions, we arrive easily at the subsystems (11-15). Solving the system (1) is equivalent to solving subsystems (11-15).
\\\\
\textbf{Proposition 3.1.} The subsystem (11) is a regular type system and its solution is given by, see \cite{2}, \cite{4}, \cite{12}, \cite{24}, \cite{25}, \cite{26}, \cite{27}
\begin{equation}
Z^p_k=J_p^{k-k_0}Z^p_{k_0}.
\end{equation}
\textbf{Proposition 3.2.} The subsystem (12) is a singular type and its solution is given by, see \cite{4}, \cite{5}, \cite{6}, \cite{7}, \cite{10}, \cite{16}
\begin{equation}
Z^q_k=0_{q,1}.
\end{equation}
\textbf{Proposition 3.3.} The subsystem (13) has infinite solutions and can be taken arbitrary
\begin{equation}
Z_k^\epsilon=C_{k,1}.
\end{equation}
\textbf{Proof.} If we set 
\begin{equation}
Z_k^\epsilon=\left[\begin{array}{c} Z_k^{\epsilon_{g+1}}\\Z_k^{\epsilon_{g+2}}\\\vdots\\Z_k^{\epsilon_d}\end{array}\right],
\end{equation}
by using (4), (5) and (20), the system (13) can be written as:
\begin{equation}
blockdiag\left\{L_{\epsilon_{g+1}}, ..., L_{\epsilon_d}\right\}\left[\begin{array}{c} Z_{k+1}^{\epsilon_{g+1}}\\Z_{k+1}^{\epsilon_{g+2}}\\\vdots\\Z_{k+1}^{\epsilon_d}\end{array}\right]=blockdiag\left\{\bar L_{\epsilon_{g+1}}, ..., \bar L_{\epsilon_d}\right\}\left[\begin{array}{c} Z_k^{\epsilon_{g+1}}\\Z_k^{\epsilon_{g+2}}\\\vdots\\Z_k^{\epsilon_d}\end{array}\right].
\end{equation}
Then for the non-zero blocks, a typical equation from (21) can be written as 
\[
\begin{array}{cc} L_{\epsilon_i} Z_{k+1}^{\epsilon_i}=\bar L_{\epsilon_i} Z_k^{\epsilon_i}, &i =g+1, g+2, ..., d, \end{array}
\]
or, equivalently,
\[
\left[\begin{array}{ccc} I_{\epsilon_i} & \vdots & 0_{{\epsilon_i}, 1}\end{array}\right]Z_{k+1}^{\epsilon_i}=\left[
\begin{array}{ccc} 0_{{\epsilon_i}, 1} & \vdots & I_{\epsilon_i}
\end{array}
\right]Z_k^{\epsilon_i},
\]
or, equivalently,
\[
\left[\begin{array}{ccccc} 1 &0&\ldots&0&0\\0&1&\ldots&0&0\\ \vdots &\vdots &\ldots&\vdots&\vdots\\0&0&\ldots&1&0\end{array}\right]\left[\begin{array}{c} z_{k+1}^{{\epsilon_i},1}\\ z_{k+1}^{{\epsilon_i},2}\\ \vdots \\ z_{k+1}^{{\epsilon_i},{\epsilon_i}} \\ z_{k+1}^{{\epsilon_i},{\epsilon_i}+1}\end{array}\right]=\left[\begin{array}{ccccc} 0 &1&\ldots&0&0\\0&0&\ldots&0&0\\ \vdots &\vdots &\ldots&\vdots&\vdots\\0&0&\ldots&0&1\end{array}\right]\left[\begin{array}{c} z_k^{{\epsilon_i},1}\\ z_k^{{\epsilon_i},2}\\ \vdots \\ z_k^{{\epsilon_i},{\epsilon_i}} \\ z_k^{{\epsilon_i},{\epsilon_i}+1}\end{array}\right],
\]
or, equivalently,
\begin{equation}
\begin{array}{ccccc}  z_{k+1}^{{\epsilon_i},1}=z_k^{{\epsilon_i},2},\\ z_{k+1}^{{\epsilon_i},2}=z_k^{{\epsilon_i},3},\\\vdots,\\ z_{k+1}^{{\epsilon_i},{\epsilon_i}}=z_{k+1}^{{\epsilon_i},{\epsilon_i}+1}.\end{array}
\end{equation}
The system (22) is a regular type system of difference equations with ${\epsilon_i}$ equations and ${\epsilon_i}+1$ unknowns. It is clear from the above analysis that in every one of the $d-g$ subsystems one of the coordinates of the solution has to be arbitrary by assigned total. The solution of the system can be assigned arbitrary
\[
Z_k^\epsilon=C_{k,1}.
\]
\textbf{Proposition 3.4.} The solution of the system (14) is unique and it is the zero solution
\begin{equation}
Z_k^\zeta=0_{t-h, 1}.
\end{equation}
\textbf{Proof.} If we set 
\[
Z_k^\zeta=\left[\begin{array}{c} Z_k^{\zeta_{h+1}}\\Z_k^{\zeta_{h+2}}\\\vdots\\Z_k^{\zeta_t}\end{array}\right],
\]
by using (6), (7) and the above expression, the system (14) can be written as
\[
blockdiag\left\{L_{\zeta_{h+1}}, ..., L_{\zeta_t}\right\}\left[\begin{array}{c} Z_{k+1}^{\zeta_{h+1}}\\Z_{k+1}^{\zeta_{h+2}}\\\vdots\\Z_{k+1}^{\zeta_t}\end{array}\right]=blockdiag\left\{\bar L_{\zeta_{h+1}}, ..., \bar L_{\zeta_t}\right\}\left[\begin{array}{c} Z_k^{\zeta_{h+1}}\\Z_k^{\zeta_{h+2}}\\\vdots\\Z_k^{\zeta_t}\end{array}\right].
\]
Then for the non-zero blocks we have:
\[
\begin{array}{ccc} L_{\zeta_j} Z_{k+1}^{\zeta_j}=\bar L_{\zeta_j} Z_k^{\zeta_j} & , & j=h+1, h+2, ..., t \end{array},
\]
or, equivalently,
\[
\left[\begin{array}{c} I_{\zeta_j} \\ \cdots \\  0_{1, {\zeta_j}}\end{array}\right]Z_{k+1}^{\zeta_j}=\left[
\begin{array}{c} 0_{1,{\zeta_j}} \\ \cdots \\ I_{\zeta_j}
\end{array}
\right]Z_k^{\zeta_j},
\]
or, equivalently,
\[
\left[\begin{array}{cccc} 1&0&\ldots&0\\0&1&\ldots&0\\ \vdots &\vdots &\ldots&\vdots\\0&0&\ldots&1\\0&0&\ldots&0\end{array}\right]\left[\begin{array}{c} z_{k+1}^{{\zeta_j},1}\\ z_{k+1}^{{\zeta_j},2}\\ \vdots \\ z_{k+1}^{{\zeta_j},{\zeta_j}}\end{array}\right]=\left[\begin{array}{cccc} 0 &0&\ldots&0\\1&0&\ldots&0\\ \vdots &\vdots &\ldots&\vdots\\0&0&\ldots&0\\0&0&\ldots&1\end{array}\right]\left[\begin{array}{c} z_k^{{\zeta_j},1}\\ z_k^{{\zeta_j},2}\\ \vdots \\ z_k^{{\zeta_j},{\zeta_j}}\end{array}\right],
\]
or, equivalently,
\begin{equation}
\begin{array}{c}  z_{k+1}^{{\zeta_j},1}=0,\\ z_{k+1}^{{\zeta_j},2}=z_k^{{\zeta_j},1},\\\vdots,\\ z_{k+1}^{{\zeta_j},{\zeta_j}}=z_k^{{\zeta_j},{\zeta_j}-1,{\zeta_j}-1},\\0=z_k^{{\zeta_j},{\zeta_j}}.\end{array}
\end{equation}
Because of the structure of $(L_{\zeta_j}, \bar L_{\zeta_j})$ blocks, it is readily shown that the only solution of (24) is the zero solution
\[
Z_k^\zeta=0_{t-h, 1}.
\]
\textbf{Proposition 3.5.} The subsystem (15) has an infinite number of solutions that can be taken arbitrary
\begin{equation}
Z_k^g=C_{k,2}.
\end{equation}
We can now state the following Theorem:
\\\\
\textbf{Theorem 3.1.} Consider the system (1), with known boundary conditions (2) and a singular matrix pencil $sF-G$. Then: 
\begin{enumerate}
\item There exists at least one solution if and only if 
\begin{enumerate}[(a)]
\item the c.m.i. of the pencil are zero, i.e. $dim\mathcal{N}_r(sF-G)=0$;
\item one of the following is satisfied
\begin{enumerate}[(i)]
\item $D\in colspan[AQ_p+BQ_pJ_p^{k_N-k_0}]$,
\item $n=p$, $D\notin colspan[AQ_p+BQ_pJ_p^{k_N-k_0}]$, $rank[AQ_p+BQ_pJ_p^{k_N-k_0}]=p$.
\end{enumerate}
Where $J_p$, $Q_p$, matrices as defined in (3), (10).
\end{enumerate}
\item If there exists a solution, then the solution is unique if and only if one of the following is satisfied
\begin{enumerate}[(i)]
\item $n=p$,$rank[AQ_p+BQ_pJ_p^{k_N-k_0}]=p$,
\item $n>p$, $D\in colspan[AQ_p+BQ_pJ_p^{k_N-k_0}]$, $rank[AQ_p+BQ_pJ_p^{k_N-k_0}]=p$.
\end{enumerate}
In this case the formula of the unique solution is given by
\begin{equation}
Y_k=Q_pJ_p^{k-k_0}Z_{k_0}^p,
\end{equation}
where $Z_{k_0}^p$ is the unique solution of the linear system
\begin{equation}
[AQ_p+BQ_pJ_p^{k_N-k_0}]Z_{k_0}^p=D.
\end{equation}
\end{enumerate}
\textbf{Proof.} First we consider that the system has non-zero c.m.i and non-zero r.m.i. By using the transformation (16), and the solutions (17-19), (23) and (25), of the subsystems (11-15) respectively, we get
\[
Z_k=\left[\begin{array}{c}
     Z_k^p  \\
     Z_k^q\\Z_k^\epsilon
     \\Z_k^\delta
     \\Z_k^g
     \end{array}\right] =\left[\begin{array}{c}
     J_p^{k-k_0}Z_{k_0}^p  \\
     0_{q, 1}\\C_{k,1}
     \\0_{t-h, 1}
     \\C_{k,2}
     \end{array}\right].
		\]
 Then
\[
     Y_k = QZ_k =
     \left[\begin{array}{ccccc}Q_p & Q_q &Q_\epsilon & Q_\zeta & Q_g \end{array}\right]
     \left[\begin{array}{c}
     J_p^{k-k_0}Z_{k_0}^p  \\
     0_{q, 1}\\C_{k,1}
     \\0_{t-h, 1}
     \\C_{k,2}\end{array}\right], 
    \]
		or, equivalently,
    \[
    Y_k =
    Q_pJ_p^{k-k_0}Z_{k_0}^p+Q_\epsilon C_{k,1}+Q_g C_{k,2}.
\]
Since $C_{k,1}$ and $C_{k,2}$ can be taken arbitrary, it is clear that the above given solution does not satisfy the given boundary conditions (2). Hence, in this case the BVP (1), (2) does not have any solutions. It is clear that the existence of c.m.i. is the reason that the systems (13) and consequently (15) exist. These systems as shown in Propositions 3.3 and 3.5 have always infinite solutions. Thus a necessary condition for the BVP to have solution is not to have any c.m.i. which is equal to 
\[
dim\mathcal{N}_r(sF-G)=0.
\]
In this case (3) will take the form
\[
sF_K  - G_K :=sI_p  - J_p  \oplus sH_q  - I_q \oplus sF_{\zeta}-G_{\zeta}.
\]
System (1) will then be is divided into the three subsystems (11), (12), (14) with solutions (17), (18), (23) respectively. Thus
\[
     Y_k = QZ_k =
     \left[\begin{array}{ccc}Q_p & Q_q & Q_\zeta\end{array}\right]
     \left[\begin{array}{c}
     J_p^{k-k_0}Z_{k_0}^p  \\
     0_{q, 1}
     \\0_{t-h, 1}
     \end{array}\right], 
    \]
		or, equivalently,
    \[
Y_k =Q_pJ_p^{k-k_0}Z_{k_0}^p.
\]
The above solution exists if and only if 
\[
AY_{k_0}+BY_{k_N}=D,
\]
or, equivalently,
\[
AQ_pZ_{k_0}^p+BQ_pJ_p^{k_N-k_0}Z_{k_0}^p=D,
\]
or, equivalently,
\[
[AQ_p+BQ_pJ_p^{k_N-k_0}]Z_{k_0}^p=D.
\]
The above linear system contains $n$ equations and $p$ unknowns. Hence, there exists at least one solution if and only if one of the following is satisfied
\[
D\in colspan[AQ_p+BQ_pJ_p^{k_N-k_0}],
\]
or,
\[
n=p,\quad D\notin colspan[AQ_p+BQ_pJ_p^{k_N-k_0}],\quad rank[AQ_p+BQ_pJ_p^{k_N-k_0}]=p.
\]
Furthermore, when there exist a solution for the linear system (27), it is unique if and only if one of the following is satisfied 
\[
p=n,\quad rank[AQ_p+BQ_pJ_p^{k_N-k_0}]=p,
\]
or,
\[
n>p,\quad D\in colspan[AQ_p+BQ_pJ_p^{k_N-k_0}],\quad rank[AQ_p+BQ_pJ_p^{k_N-k_0}]=p.
\]
The unique solution of the BVP (1), (2) is then given by (26). The proof is completed.

\section*{Non-consistent boundary value problem}

From Theorem 3.1, a singular boundary value problem (BVP) of type (1)-(2) can have a unique solution, infinite solutions or no solutions.  
\\\\
\textbf{Definition 3.1.} If for the system (1) with boundary conditions (2) there exists at least one solution, the BVP (1)-(2) is said to be consistent.
\\\\
We can now state the following Theorem, which provides optimal solutions for non-consistent BVPs of type (1)-(2). With $\left\|\cdot\right\|$ we will denote an induced norm, with $\left\|\cdot\right\|_2$ the euclidean norm and with $()^*$ the conjugate transpose tensor.
\\\\
\textbf{Theorem 3.2.} We consider the BVP (1)-(2) with a singular pencil, $\left\|J_p\right\|<1$ and all column minimal indices zero, i.e. $dim \mathcal{N}_r(sF-G)=0$. Let $K=AQ_p+BQ_pJ_p^{k_N-k_0}$, then for the non-consistent BVP of type (1)-(2)
\begin{enumerate}[(a)]
\item If $p< n$, $D\notin colspan K$ and $K$ is full rank, an optimal solution of the BVP (1)-(2) is given by
\begin{equation}
    \hat Y_k=Q_pJ_p^{k-k_0}(K^*K)^{-1}K^*D.
\end{equation}
\item If $D\notin colspan K$ and $K$ is rank deficient, an optimal solution of the BVP (1)-(2) is given by
\begin{equation}
    \hat Y_k=Q_pJ_p^{k-k_0}(K^*K+E^*E)^{-1}K^*D.
\end{equation}
Where $E$ is a matrix such that $K^*K+E^*E$ is invertible and $\left\|E\right\|_2=\theta$, $0<\theta<<1$.
\end{enumerate}
\textbf{Proof.} From Theorem 3.1, $Z_{k_0}^p$ is the solution of the linear system
\begin{equation}
KC=D,
\end{equation}
where $C\in\mathbb{C}^{p\times 1}$ the unknown matrix. System (30) has $n$ linear equations and $p$ unknowns. 
\\\\
For the the proof of (a), since $p< n$, $D\notin colspan K$ and $K$ is full rank, the system (30) has no solutions. Let $\hat D$ be a vector such that a vector $\hat C$ is the unique solution of the system $K\hat C=\hat D$. We want then to solve the following optimization problem
\[
\begin{array}{c}min\left\|D-\hat D\right\|_2^2\\s.t.\quad K\hat C=\hat D,\end{array}
\]
or, equivalently,
\[
min\left\|D-K\hat C\right\|_2^2.
\]
Where $\hat C$ is the optimal solution, in terms of least squares, of the linear system (30). In this case, the solution $\hat C$ is given by
\[
\hat C=(K^*K)^{-1}K^*L.
\]
This is the least squares solution and hence an optimal solution of system (1) with boundary conditions of type (2) is given by (28). 
\\\\
For the proof of (b), since $D\notin colspanK$ and $K$ is rank deficient, the system (30) has no solutions. While the matrix $K$ is rank deficient, the matrix $K^*K$ is singular and hence not invertible. Thus we can not apply the same method as in (a). In this case, we seek a solution $\hat C$ minimizing the functional
\[
D_2(\hat C)=\left\|D-K\hat C\right\|_2^2+\left\|E\hat C\right\|_2^2.
\]
Where $E$ is a matrix such that $K^*K+E^*E$ is invertible and $\left\|E\right\|_2=\theta$, $0<\theta<<1$. Expanding $D_2(\hat C)$ gives
\[
D_2(\hat C)=(D-K\hat C)^*(D-K\hat C)+(E\hat C)^*E\hat C,
\]
or, equivalently,
\[
D_2(\hat C)=D^*D-2D^*K\hat C+(\hat C)^*K^*K\hat C+(\hat C)^*E^*E\hat C,
\]
because $D^*K\hat C=(\hat C)^*K^*D$. Furthermore
\[
\frac{\partial}{\partial \hat C }D_2(\hat C)=-2K^*D+2K^*K\hat C+2E^*E\hat C.
\]
Setting the derivative to zero, $\frac{\partial}{\partial \hat C }D_2(\hat C)=0$, we get
\[
(K^*K+E^*E)\hat C=K^*D.
\]
The solution is then given by
\[
\hat C=(K^*K+E^*E)^{-1}K^*D.
\]
Hence an optimal solution in this case is given by (29). The proof is completed.
\\\\
\textbf{Remark 3.1.} In Theorem 3.2 the optimal solution of the BVP problem (1)-(2) was achieved after a perturbation to the column $D$ accordingly 
\[
min\left\|D-\hat D\right\|_2,
\]
or, equivalently,
\[
\left\|D-K(K^*K)^{-1}K^*D\right\|_2,
\]
in case of (a) and 
\[
\left\|D-K(K^*K+E^*E)^{-1}K^*D\right\|_2
\]
in case of (b). 
\\\\
For a consistent BVPs of type (1)-(2) with infinite solutions, there are different ways to obtain optimal solutions, depending on the problem that the BVP represents. In the following Proposition, we apply some techniques in order to gain the minimum solution.
\\\\
\textbf{Proposition 3.6.} We consider the BVP (1)-(2) with a singular pencil, $\left\|J_p\right\|<1$ and all column minimal indices are zero, i.e. $dim \mathcal{N}_r(sF-G)=0$. Let $K=AQ_p+BQ_pJ_p^{k_N-k_0}$. Then, for a consistent BVP of type (1)-(2) with infinite solutions:
\begin{enumerate}[(a)]
\item If $K^\dagger$ is the Moore-Penrose Pseudoinverse of $K$, then an optimal solution of the BVP (1)-(2) is given by 
\begin{equation}
    \hat Y_k=Q_pJ_p^{k-k_0}K^\dagger D.
\end{equation}
\item If $n> p$, $K$ is full rank, a minimum solution of the BVP (1)-(2) is given by
\begin{equation}
    \hat Y_k=Q_pJ_p^{k-k_0}(KK^*)^{-1}KD.
\end{equation}
\item If $D\in colspanK$ and $K$ is rank deficient, an optimal solution of the BVP (1)-(2) is given by (29).
\end{enumerate}
\textbf{Proof.} From Theorem 3.1, $Z_{k_0}^p$ is the solution of the linear system (30) which has $n$ linear equations and $p$ unknowns. 
\\\\
For the proof of (a), given an $n\times p$ matrix $K$, the Moore-Penrose Pseudoinverse of $K^\dagger$ is calculated from the singular value decomposition of $K$, see \cite{24}, \cite{25}, \cite{26}, \cite{27}. The solution of system (30) is then given by $\hat C=K^\dagger L$.
\\\\
For the proof of (b), $p> n$, $K$ is a wide matrix (more columns than rows) and full rank, i.e. the system (30) is underdetermined. In this case, it is common to seek a solution $\hat C$ (one of the infinite solutions of system (30)) with minimum norm. Hence we would like to solve the optimization problem
\[
\begin{array}{c}\quad min \quad \left\|\hat C\right\|_2^2,\\
\quad s.t.\quad D=K\hat C.
\end{array}
\]
By defining the Lagrangian
\[
D_3(\hat C,\lambda)=\left\|\hat C\right\|_2^2+\lambda^*(D-K\hat C)
\]
and by taking the derivatives of the Lagrangian we get
\[
\frac{\partial}{\partial \hat C }D_3(\hat C,\lambda)=2\hat C- K^*\lambda
\]
and
\[
\frac{\partial}{\partial \lambda }D_3(\hat C,\lambda)=D-K\hat C.
\]
Setting the derivatives to zero, $\frac{\partial}{\partial \hat C }D_3(\hat C,\lambda)=0$ and $\frac{\partial}{\partial \lambda }D_3(\hat C,\lambda)=0$, we get
\[
\hat C=\frac{1}{2} K^*\lambda
\]
and
\[
D=K\hat C,
\]
or, equivalently,
\[
D=\frac{1}{2}KK^*\lambda.
\]
Since $rank K= n$, the matrix $KK^*$ is invertible and thus
\[
\lambda=2(KK^*)^{-1}D.
\]
The solution is then given by
\[
\hat C=K^*(KK^*)^{-1}D.
\]
Hence an optimal solution in this case is given by (32).
\\\\
For the proof of (c), since $D\in colspanK$ and $K$ is rank deficient, the system (30) has infinite solutions. While matrix $K$ is rank deficient, the matrix $KK^*$ is singular and hence not invertible. In this case, we seek a solution $\hat C$ by minimizing the functional
\[
D_4(\hat C)=\left\|D-K\hat C\right\|_2^2+\left\|E\hat C\right\|_2^2,
\]
which is given by
\[
\hat C=(K^*K+E^*E)^{-1}K^*D
\]
and hence an optimal solution in this case is given by (29). The proof is completed.

\section{Numerical examples}

In \cite{2} we get a description of the Leontief model. The Leontief model is a dynamic model of a multisector economy. One of its versions is constructed as follows. Suppose the economy is divided into $m$ sectors. Let $Y_k$ have $i^{th}$ component the output in the $k^{th}$ time period of sector $i$. Let time $k$ be finite, $k_0\leq k\leq k_N$, i.e. we want our investigation to be between the time period $k_0$ and $k_N$. Let $f_{ij}$ be the amount of commodity $i$ that sector $j$ must have to produce one unit of commodity $j$. Let $f_{ij}$ be elements of the matrix $F$. Let $m_{ij}$ be the proportion of commodity $j$ that gets transferred to commodity $i$ in the $k_{th}$ time period and $m_{ij}$ the elements of the matrix $M$ (also called the Leontief input-out matrix or the matrix of flow coefficients). The Leontief model then says that amount of commodity $i$ is equal to ammount of commodity $i$ inputed by all sectors plus ammount needed for production of the next output $Y_{k+1}$. In matrix form
\[
Y_k=MY_k+F(Y_{k+1}-Y_k).
\]
Where $M$ is called the Leontief input-output coefficient matrix or the matrix of flow coefficients, $F$ is the capital coefficient matrix, $Y_k$ is the vector of output levels in period $k$. The first term on the right-hand side of the above notation, $MY_k$, denotes intermediate demand for goods by industries; whereas the second term, $F(Y_{k+1}-Y_k)$, reflects the distribution of inputs to investment. The capital coefficient matrix $F$ is usually singular because only a small number of sectors ordinarily contribute to capital formation. For $G=I_m-M+F$ we arrive at the system (1).
\\\\    
\textit{Example 1}\\\\
Consider the singular system (1) and let 
\[
F=\left[\begin{array}{ccccccc} 2&1&1&0&0&0&0\\1&3&1&1&0&0&0\\1&1&2&1&0&0&0\\0&1&1&1&0&0&0\\0&0&0&0&0&0&0\\0&0&0&0&1&0&0\\0&1&0&0&0&0&1\end{array}\right],
G=\left[\begin{array}{ccccccc} 1&1&1&0&0&0&1\\0&3&2&2&0&1&1\\1&2&3&2&0&0&0\\0&2&2&2&0&0&0\\0&0&0&0&1&0&0\\0&0&0&0&0&0&0\\0&0&0&0&0&1&0\end{array}\right].
\]
Then det$[sF-G]$=0. The invariants of the pencil are,  $s-2$, $s-1$ (finite elementary divisors), $\epsilon_1$=0, $\epsilon_2$=2 the (column minimal indices) and 
$\zeta_1$=0, $\zeta_2$=1 are the (row minimal indices). From Theorem 3.1 there does not exist a solution since the c.m.i. are non-zero.
\\\\
\textit{Example 2}\\\\
Consider the singular system (1) for $k=0,1,...,100$ and let 
\[
F=\left[\begin{array}{ccccc} 1&1&1&1&1\\0&1&1&0&1\\1&1&1&1&1\\0&1&1&0&1\\1&0&1&0&0\\0&0&1&1&1\end{array}\right],G=\left[\begin{array}{ccccc} 1&2&2&1&2\\0&2&2&0&2\\1&2&2&2&3\\0&2&3&1&3\\0&0&0&0&0\\1&0&1&0&0\end{array}\right].
\]
The matrices $F$, $G$ are non square. Thus the matrix pencil $sF-G$ is singular with invariants, the finite elementary divisors $s-2$, $s-1$, 
an infinite elementary divisor of degree 1 and the row minimal indices $\zeta_1$=0, $\zeta_2$=1. From Theorem 3.1 there exist non-singular matrices
\[
P=\left[\begin{array}{cccccc}1&-1&0&0&0&0\\0&1&0&0&0&0\\-1&0&1&0&0&0\\0&0&0&0&1&0\\0&0&0&0&0&1\\0&-1&0&1&0&0\end{array}\right],
Q=\left[\begin{array}{ccccc}0&0&1&1&-1\\1&1&-1&-1&0\\0&0&-1&0&1\\1&0&-1&-1&1\\-1&0&2&1&-1\end{array}\right],
\]
such that $PFG=F_K$ and $PGQ=G_K$. Where
\[
F_K = \left[\begin{array}{ccccc}1&0&0&0&0\\0&1&0&0&0\\0&0&0&0&0\\0&0&0&1&0\\0&0&0&0&1\\0&0&0&0&0\end{array}\right],
G_K=\left[\begin{array}{ccccc}1&0&0&0&0\\0&2&0&0&0\\0&0&1&0&0\\0&0&0&0&0\\0&0&0&1&0\\0&0&0&0&1\end{array}\right].
\]
Hence
\[
Q_p=\left[\begin{array}{cc}0&0\\1&1\\0&0\\1&0\\-1&0\end{array}\right],J_p=\left[\begin{array}{cc}1&0\\0&2\end{array}\right].
\]
Let the system (1) have boundary conditions of type (2) with 
\[
A=\left[\begin{array}{ccccc}1&0&0&0&0\\1&1&0&0&0\\1&0&0&0&0\\1&0&0&2&1\\1&0&0&2&1\end{array}\right],
B=\left[\begin{array}{ccccc}1&0&0&1&1\\2&0&2&2&2\\1&0&3&1&1\\3&0&1&0&0\\0&0&0&1&0\\4&0&2&1&1\end{array}\right], D=\left[\begin{array}{c}0\\-1\\0\\1\\-1\end{array}\right].
\]
Then 
\[
AQ_p+BQ_pJ_p^{k_N-k_0}=\left[\begin{array}{cc}0&0\\1&1\\0&0\\1&0\\-1&0\end{array}\right].
\]
Consequently, since $n=5>2=p$ with 
\[
rank[AQ_p+BQ_pJ_p^{k_N-k_0}]=rank\left[\begin{array}{cc}0&0\\1&1\\0&0\\1&0\\-1&0\end{array}\right]=2
\]
and
\[
D\in colspan\left[\begin{array}{cc}0&0\\1&1\\0&0\\1&0\\-1&0\end{array}\right],
\]
from Theorem 3.1 the solution of the BVP (1), (2) is
\[
Y_k=\left[\begin{array}{cc}0&0\\1&1\\0&0\\1&0\\-1&0\end{array}\right]\left[
     \begin{array}{cc} 1&0\\0&2^k\end{array}\right]Z_0^p
\]
and by calculating $Z_0^p$ we get
\[
Z_0^p=\left[\begin{array}{c}
 1\\-2
\end{array}\right].
\]
Hence, the unique solution of the BVP (1), (2) will be 
\[
Y_k=\left[\begin{array}{c} 0\\1-2^{k+1}\\0\\1\\-1\end{array}\right].
\]
\textit{Example 3}\\\\
Consider the BVP (1), (2) for $k=0,1,...,100$ and let $F,G,B$ be defined as in Example 2. Furthermore let
\[
A=\left[\begin{array}{ccccc}1&0&0&0&0\\1&1&0&0&0\\1&0&0&0&0\\1&1&0&1&1\\1&-1&0&2&2\end{array}\right],D=\left[\begin{array}{c}  0\\1\\0\\1\\-1\end{array}\right].
\]
The pencil $sF-G$ has the same invariants as in example 2, i.e. the finite elementary divisors $s-2$, $s-1$, an infinite elementary divisor of degree 1 and the row minimal indices $\zeta_1$=0, $\zeta_2$=1. Then 
\[
AQ_p+BQ_pJ_p^{k_N-k_0}=\left[\begin{array}{cc}0&0\\1&1\\0&0\\1&1\\-1&-1\end{array}\right].
\]
Consequently, since there are no column minimal indices and
\[
D\in colspanAQ_p+BQ_pJ_p^{k_N-k_0},
\]
from Theorem 3.1 there exist solutions for the BVP (1), (2), but because $n=5>2=p$ with 
\[
rank[AQ_p+BQ_pJ_p^{k_N-k_0}]=rank\left[\begin{array}{cc}0&0\\1&1\\0&0\\1&1\\-1&-1\end{array}\right]=1<2=p,
\]
the solution is not unique. Hence, 
\[
Y_k=\left[\begin{array}{cc}0&0\\1&1\\0&0\\1&0\\-1&0\end{array}\right]\left[
     \begin{array}{cc} 1&0\\0&2^k\end{array}\right]Z_0^p
\]
and by calculating $Z_0^p$ we get
\[
Z_0^p=\left[\begin{array}{c}
 1-c\\c
\end{array}\right].
\]
Thus, the solutions of the BVP (1), (2) will be given by
\[
\begin{array}{cc}Y_k=\left[\begin{array}{c} 0\\1-c-2^{k+1}c\\0\\1-c\\-1+c\end{array}\right],&\forall c\in \mathbb{C}\end{array}.
\]
From Proposition 3.6 an optimal solution is given by (29)
\[
\hat Y_k=Q_pJ_p^{k-k_0}([AQ_p+BQ_pJ_p^{k_N-k_0}]^*[AQ_p+BQ_pJ_p^{k_N-k_0}]+E^*E)^{-1}[AQ_p+BQ_pJ_p^{k_N-k_0}]^*D.
\]
Where $E=\left[\begin{array}{cc}0&0\\1&\theta\end{array}\right]$ is a matrix such that $[AQ_p+BQ_pJ_p^{k_N-k_0}]^*[AQ_p+BQ_pJ_p^{k_N-k_0}]+E^*E$ is invertible and $\left\|E\right\|_2=\theta$, $0<\theta<<1$. Hence
\[
\hat Y_k=\left[\begin{array}{cc}0&0\\1&1\\0&0\\1&0\\-1&0\end{array}\right]\left[
     \begin{array}{cc} 1&0\\0&2^k\end{array}\right](\left[\begin{array}{cc}0&0\\1&1\\0&0\\1&1\\-1&-1\end{array}\right]
^*\left[\begin{array}{cc}0&0\\1&1\\0&0\\1&1\\-1&-1\end{array}\right]
+E^*E)^{-1}\left[\begin{array}{cc}0&0\\1&1\\0&0\\1&1\\-1&-1\end{array}\right]
^*\left[\begin{array}{c}  0\\1\\0\\1\\-1\end{array}\right],
\]
or, equivalently,
\[
\hat Y_k=\left[\begin{array}{c}  0\\1\\0\\1\\-1\end{array}\right].
\]
\section*{Conclusions}

In this article, we study a BVP of a class of linear singular discrete time systems whose coefficients are either non-square constant matrices or square matrices with a pencil which has an identically zero determinant. Actually, the results of the paper apply to general non-square pencils and generalize various results regarded the literature which mainly are dealing with square and non-singular systems. Another important characteristic of the general case considered here is that the uniqueness of solutions is not automatically satisfied. This is very important for many applications for which the model is significant only for certain range of its parameters. In these cases a careful interpretation of results or even a redesign of the system maybe needed. By taking into consideration that the relevant pencil of system (1) is singular, we decompose the linear singular matrix difference equation into five sub-systems. Afterwards, we provide necessary and sufficient conditions for existence and uniqueness of solutions. More analytically we analyze the conditions under which the system has unique and infinite solutions. Explicit and easily testable conditions are derived for which the system has a unique solution. Furthermore, we provide a formula for the case of the unique solution and study optimal solutions for the cases of no solutions and infinite many solutions. Finally, we apply our results to a reformulated Leontief model.

\section*{Acknowledgement}

I. Dassios is funded by Science Foundation Ireland (award 09/SRC/E1780).

\end{document}